\documentclass[a4paper,letterpaper,leqno]{amsart}

\usepackage[english]{babel}
\usepackage{url}
\usepackage[all,cmtip]{xy}
\usepackage{color}
\usepackage{hyperref}
\usepackage{enumerate}
\usepackage{array}
\usepackage{nicefrac}
\usepackage[section]{placeins}
\usepackage[linguistics]{forest}
\usepackage{amsmath,amssymb,enumerate,amsthm}
\usepackage{cite}
\usepackage{latexsym}
\usepackage{setspace}
\usepackage{mathtools}
\usepackage[labelfont=bf]{caption}
\usepackage{longtable}

\theoremstyle{plain}
\newtheorem{theorem}{Theorem}[section]
\newtheorem{nonumtheorem}{Theorem}

\newtheorem{lemma}[theorem]{Lemma}

\theoremstyle{remark}
\newtheorem{remark}{Remark}[section]

\numberwithin{equation}{section}

\def\A{\mathcal{A}}
\def\F{\mathcal{F}}
\def\H{\mathcal{H}}
\def\N{\mathbb{N}}

\def\R{\mathbb{R}}
\def\Z{\mathbb{Z}}

\def\SL{\mathrm{SL}(2,\mathbb{Z})}
\def\GL{\mathrm{GL}(2,\mathbb{Z})}

\def\Im{\mathrm{Im}}
\def\Re{\mathrm{Re}}

\def\={\;=\;}
\def\.={\;\dot{=}\;}
\def\l{\ell}

\begin{document}

\title[Cycle integrals of the $j$-function on Markov geodesics]{Asymptotic bounds for  cycle integrals of the $j$-function on Markov geodesics}
\author{Paloma Bengoechea}
\address{ETH Zurich \\ Department of Mathematics\\ Ramistrasse 101\\8092 Zurich\\Switzerland}
\email{paloma.bengoechea@math.ethz.ch}
\thanks{Bengoechea's research is supported by SNF grant 173976.}

\maketitle

\begin{abstract}
	We give asymptotic
	 upper and lower bounds for the real and imaginary parts of  cycle integrals of  the classical modular $j$-function 
	 along  geodesics that correspond to Markov irrationalities.
	
	\end{abstract}

\section{Introduction}
Let $w$ be a real quadratic irrationality, so a root of an equation  
$$
ax^2+bx+c=0\qquad(a,b,c\in\Z,\quad (a,b,c)=1)
$$
with positive non square discriminant $D=b^2-4ac$. Let  $Q$ be the quadratic form $[a,b,c]$. 
The geodesic $S_Q$ in the upper half plane $\H$ joining $w$ and its Galois conjugate is given by the equation
$$
a|z|^2+b\, \mathrm{Re}\, z+c=0 \qquad (z\in\mathcal{H}).
$$
We orientate $S_Q$  counterclockwise if $a>0$ and clockwise if $a<0$. 
There is an infinite cyclic group $\Gamma_Q$ in $\SL$, corresponding to the group of totally positive units in $\mathbb{Q}(\sqrt{D})$, that preserves $Q$ and hence $S_Q$. The smallest positive unit  in $\mathbb{Q}(\sqrt{D})$ is given by $\varepsilon=(t+u\sqrt{D})/2$, where $(t,u)$ is the smallest positive integral solution to Pell's equation $t^2-Du^2=4$.
 The induced geodesic $C_Q=\Gamma_Q\backslash S_Q$  on the modular surface is closed, primitive, positively-oriented with length
 $$
\mathrm{length}(C_Q)=\int_{ C_Q} \frac{\sqrt{D}}{Q(z,1)}dz =2\log \varepsilon.
$$
Here $\frac{\sqrt{D}}{Q(z,1)}dz$ is the hyperbolic arc length on $S_Q$.
As usual, we denote by $j(z)$ the modular function for $\SL$  introduced by 
 Klein as an invariant of elliptic curves. 
The group  $\Gamma_Q$   preserves  $j(z)\frac{\sqrt{D}}{Q(z,1)}dz$, and so  we can consider the integral
$$
 \int_{C_Q} j(z) \frac{\sqrt{D}}{Q(z,1)}dz
$$
and define  the \textit{`value'} (also called \textit{cycle integral}) of $j$ at $w$,  as the complex number:
\begin{equation}	\label{defcycle}
j(w)=j(Q):=\dfrac{1}{2\log\varepsilon} \int_{C_Q} j(z)\frac{\sqrt{D}}{Q(z,1)}dz.
\end{equation}
Note that the integral above  is  $\SL$-invariant.

Cycle integrals  have been related to mock modular forms \cite{DIT1}, to modular knots \cite{DIT-link} and to class numbers of real quadratic
fields \cite{DIT-hurw}.
Moreover, cycle integrals of the Klein invariant $j$ share several analogies with singular moduli (the values of the $j$-function at imaginary quadratic irrationalities) when 
both are  gathered in `traces' (see \cite{DIT1}, \cite{DIT2}, \cite{DFI},
\cite{masri}). 
By analogy to the traces of singular moduli (the values of the $j$-function at imaginary quadratic irrationalities) 
we define $\mbox{Tr}_D j:=\sum j(w_Q)$ where the sum is  over $\SL$-classes of indefinite binary quadratic forms $Q$ of fundamental discriminant $D>0$.
 The asymptotic distribution of the traces was studied in \cite{Duke} and \cite{DFI},\cite{masri} for negative and positive discriminants respectively.
As fundamental discriminants $D\rightarrow+\infty$, it was shown in \cite{DFI},\cite{masri} that
\begin{equation}\label{tr} \frac{\mbox{Tr}_ D (j)}{\mbox{Tr}_D(1)}\rightarrow 720.
\end{equation} 
Here $720$ is an `average' value of the $j$ function.
The individual values remain very much unknown.
 In \cite{BI2} we gave the first  bounds  for the real parts of the values $j(w)$  (in fact $j$ could be replaced  by  any modular function $f$ which is real valued on the geodesic arc $\left\{e^{i\theta} : \pi/3 \leq \theta \leq 2\pi/3\right\}$).  We showed that 
 $\Re(j(w))\leq 744$ (we recall that $744$ is the constant term in the Fourier expansion of $j$). For quadratic irrationalities $w$ satisfying a `quite strong' diophantine condition, we  proved that  $\Re(j(w))\geq j(\frac{1+\sqrt{5}}{2})\approx 706.3248$. These bounds are optimal and were conjectured earlier by Kaneko in \cite{Ka} based on numerical evidence. The imaginary parts of $j(w)$ are conjectured to lie in $(-1,1)$, but to our knowledge nothing is known yet. 
 
 The values $j(w)$ are particularly interesting at Markov irrationalities $w$. Markov irrationalities are important in diophantine approximation  and in the search of 
  positive minima of indefinite binary quadratic forms. They are structured on a tree called the Markov-Hurwitz tree and there is a very rich theory attached to them in the interplay of number theory, diophantine approximation, hyperbolic geometry, dynamics and graph theory. 
Kaneko published in \cite{Ka} the first numerical data on values $j(w)$ at Markov irrationalities, together with some conjectures based on his data. Namely, he conjectured that each $j(w)$ is between the values of $j$ at two Markov irrationalities that lie  above $w$ on the tree.  We call these two irrationalities the predecessors of $w$ and we call call this property `interlacing property'. In \cite{BI1} we proved that for every branch of the tree, the interlacing property holds after some level that depends on the branch. One can think of that result as a `local' asymptotic interlacing property (an asymptotic property that holds for each branch).

In this paper we prove the `global' interlacing property, namely:
\begin{nonumtheorem}\label{th1}
	 Let $w=w_n$ be a Markov irrationality and  $n$ be the  level of $w$ in the Markov-Hurwitz tree. There is a unique path on the Markov-Hurwitz tree that ends with $w_n$; let $w_{n-1}$ be the Markov irrationality above $w_n$ on the path, and  $w_{n'-1}$ be the Markov irrationality above the last irrationality  where the path turns. There exists an  integer $n_0$ such that, if $n> n_0$, then
	$j(w_n)$ lies between $j(w_{n-1})$ and $j(w_{n'-1})$.
\end{nonumtheorem}
The irrationalities $w_{n-1}$ and $w_{n'-1}$ are the `predecessors' of $w_n$ in Theorem \ref{th1}.
The interlacing property is expected to hold for $n\geq 2$, hence we expect the bounds:
\begin{align*}
	706.32481\approx j((1+\sqrt{5})/2)\leq  &\Re(j(w))\leq j(\sqrt{2})\approx 709.8929\\
	-0.26703\ldots\leq &\Im(j(w)) \leq 0.26703\ldots.
\end{align*}
By making explicit most of the estimations in Theorem \ref{th1}, we obtain: 
\begin{nonumtheorem}\label{th2}
	Let $w$ be a Markov irrationality and  $n$ be the  level of $w$ in the Markov-Hurwitz tree. We have that
	\begin{align*}
		681.50081\leq \Re(j(w))&\leq 	742.03641\qquad as\ n\rightarrow\infty\\
	-0.93637	\leq \Im(j(w))&\leq 0.67396	 \qquad as\ n\rightarrow\infty.
\end{align*}
\end{nonumtheorem}
Although the bounds in Theorem \ref{th2} are not optimal, they are the first lower bounds for the real parts of the values $j(w)$ at Markov irrationalities and the first bounds for the imaginary parts that we have until now.

 One reason why we obtained only a `local'  interlacing property in \cite{BI1} is  that we lacked of a convenient order on the tree.   A new key approach here is to order the Markov-Hurwitz tree by Farey fractions (which can
  probably be more exploited in the future).  We also exploit more  accurately the diophantine properties of the Markov-Hurwitz tree as well as its relations with the Farey tree, and use an asymptotic formula from \cite{ZMar} on Markov numbers.
The  arguments used in this paper apply to any modular function $f$ which is real valued on the geodesic arc $\left\{e^{i\theta} : \pi/3 \leq \theta \leq 2\pi/3\right\}$. The values \eqref{average}, \eqref{compreal}, \eqref{compimag} will change when we replace $j$ by $f$, and so will the subsequent calculations.

The paper is organized as follows. In the next two sections,  we quickly review the main key points of Markov's theory and other related facts that will be useful. Concretely, in section 2, we introduce Markov numbers and Markov irrationalities. In section 3 we introduce the Markov-Hurwitz and the Farey trees as well as  some of their properties and interrelations. In section 4 we give an asymptotic formula for $\log\varepsilon$ in terms of the Farey tree by using Zagier's  asymptotic formula for Markov numbers in \cite{ZMar}. 
In section 5 we write the values $j(w)$ in terms of cycles of certain reduced binary quadratic forms, usually known as `simple forms' after Zagier. Simple forms can be written themselves in terms of cycles of quadratic irrationalities produced by a certain continued fraction algorithm. We work with the cycles of continued fraction expansions in sections 6 and 7 to obtain a `local' formula  for $(2\log\varepsilon) j(w)$ that depends on two neighbours of $w$ on the Markov-Hurwitz tree (the two predecessors of $w$). In section 8 we deduce a `global' formula for $(2\log\varepsilon) j(w)$ from the local formula obtained  in section 7. We prove Theorem \ref{th1} in section 9 using the global formula from section 8 and  the asymptotic formula for $\log\varepsilon$ given in section 4. We prove Theorem \ref{th2} in section 10. We give some numerical data in the appendix, in section 11.

\section{Markov's theory}

Markov's work (1880, \cite{mar1}, \cite{mar2}) establishes very beautiful connections between  positive integral minima 
 for indefinite binary  quadratic forms  and the Lagrange-Hurwitz problem in diophantine approximation. 
The Lagrange-Hurwitz problem consists in describing the Lagrange constants
$$
L(x)=(\liminf_{q\rightarrow\infty} q\|qx\|),
$$
 where $x$ runs through the real numbers and $\|\cdot\|$ denotes the distance to a closest integer.
The quantity $L(x)$  provides a `measure' of how well $x$ can be approximated by  the rationals. 
 For almost all $x\in\R$ we have $L(x)=0$, and when $L(x) >0$ we call $x$ \emph{badly approximable}. 
 A well known theorem of Hurwitz states that $L(x)<1/\sqrt{5}$.

The  \textit{Lagrange spectrum} $L:=\left\{L(x)^{-1}\right\}_{x\in\R}\subseteq[\sqrt{5},\infty]$
is  structured in three parts: $\L\cap[\sqrt{5},3)$ is discrete with $3$ as the only accumulation point, 
$\L\cap(3,F]$ is fractal, where $F\approx 4.528$ is Freiman's constant, and $\L\cap(F,\infty]$ is continuous.

Markov discovered that on the interval $(0,3)$ the Lagrange spectrum and the spectrum $M$ of integral minima 
 $$
 M(Q)=\dfrac{\sqrt{D}}{\underset{\substack{(x,y)\in\Z^2\\(x,y)\neq(0,0)}}{\mathrm{\inf}} |Q(x,y)|}
 $$
of  binary quadratic forms
$$
Q(x,y)=Ax^2+Bxy+Cy^2 \qquad(A,B,C\in\R)
$$
with positive discriminant $D=B^2-4AC$, are the same.

 The unifying thread  is 
 the diophantine equation
 \begin{equation}\label{Deq}
  a^2+b^2+c^2=3abc\qquad (a,b,c\in \N),
 \end{equation}
and its integer solutions.
  The  integer solutions $(a,b,c)$ are obtained by starting with $(1,1,1)$, $(1,1,2)$,
   $(2,1,5)$ and then proceeding recursively going from $(a,b,c)$ to the new triples obtained by Vieta involutions $(c,b,3bc-a)$ and $(a,c,3ac-b)$.
    \emph{Markov numbers}  are the greatest coordinates of each solution $(a,b,c)$. They form the  \emph{Markov sequence}
      \begin{equation}\label{markovseq}
      	\left\{c_i\right\}_{i=1}^\infty=\left\{1,2,5,13,29,34,89,169,194,\ldots\right\}.
      \end{equation}
  (We count multiplicities if the unicity conjecture is false.)
To each   Markov number $c$,  a  \emph{Markov irrationality}
$$
w=\dfrac{3c-2k+\sqrt{9c^2-4}}{2c}
$$
is associated, 
where $k$ is an integer that satisfies $a k\equiv b\pmod{c}$, $0\leq k<c$,  and $(a,b,c)$ is a solution to \eqref{Deq} with $c$ maximal. 
The quadratic $w$ is a root of the \textit{Markov form}
$$
Q(x,y)=c x^2 + (3c - 2k)xy + (\ell - 3k)y^2,
$$
where $\ell=\frac{k^2+1}{c}\in\Z$. 
Then
$$
M(Q)=L(w)^{-1}=\sqrt{9-4/c^2}.
$$
Moreover,
$$
M\cap(0,3)=\L\cap(0,3)=\left\{\sqrt{9-4/c_i^2}\right\}_{i\geq 1}
$$
and any $x\in\R$ with $L(x)^{-1}\in (0,3)$ or any $Q$ with $M(Q)\in (0,3)$ is $\GL$-equivalent to a Markov quadratic or
a Markov form respectively.

  \section{Trees related to Markov's theory}
 
 \subsection{The Markov-Hurwitz tree}
 
 The solutions of the diophantine equation \eqref{Deq}  inherit a tree structure with two bifurcations from the two  Vieta involutions   described earlier:
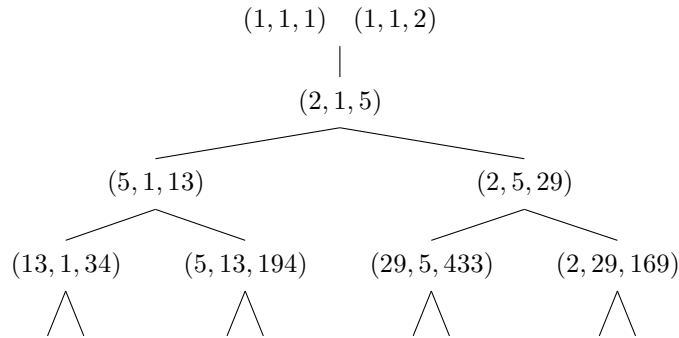
\begin{figure}[h!]
	\centering
	\begin{forest} 
[
{$(1,1,1)$\quad$(1,1,2)$}
[
{$\begin{array}{c}(2,1,5)\end{array}$} 
[
{$\begin{array}{c}(5,1,13)\end{array}$} 
[
{$\begin{array}{c}(13,1,34)\end{array}$}
[]
[]]
[
{$\begin{array}{c}(5,13,194)\end{array}$}
[]
[]
]] 
[
{$\begin{array}{c}(2,5,29)\end{array}$} 
[
{$\begin{array}{c}(29,5,433)\end{array}$}
[]
[]] 
[
{$\begin{array}{c}(2,29,169)\end{array}$}
[]
[]]]]]
\end{forest}
	\caption{The Markov tree}
\end{figure}  

 Naturally, we can consider the parallel tree of Markov irrationalities, where each vertex is a quadratic irrationality that corresponds to a Markov number.
 What is more interesting is that 
  each Markov irrationality can be constructed from two predecessors on the tree by a simple conjunction operation on the continued fraction expansions.  A formal decription of the procedure using  `+' continued fraction expansions can be found in \cite{Bom}. It will be more convenient for us to work with the `--' continued fraction since it corresponds to transformations by the modular group (whereas the `+' continued fraction corresponds to transformations by $\GL$). It was observed in \cite{BI1} that roughly the same conjunction operation works with the `--' continued fraction.
 We denote by   
$$
(a_1,a_2,a_3,\ldots)=a_1-\dfrac{1}{a_2-\dfrac{1}{a_3-\dfrac{1}{\ddots}}}
$$
 the  `--' continued fraction expansion, with $a_i\in\Z$ and $a_i \geq 2$ for $i\geq 2$. 

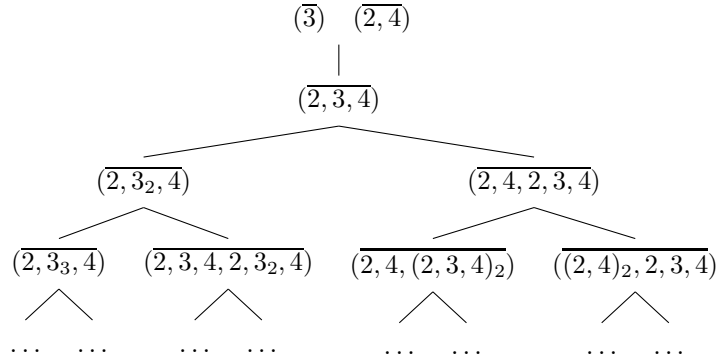
\begin{figure}[h!]
	\centering
	\begin{forest}
		[
		\quad{$(\overline{3})$\quad$(\overline{2,4})$}
		[
		{$(\overline{2,3,4})$} 
		[
		{$(\overline{2,3_2,4})$} 
		[
		{$(\overline{2,3_3,4})$}
		[$\ldots$]
		[$\ldots$]]
		[
		{$(\overline{2,3,4,2,3_2,4})$}
		[$\ldots$]
		[$\ldots$]
		]] 
		[
		{$(\overline{2,4,2,3,4})$} 
		[
		{$(\overline{2,4,(2,3,4)_2})$}
		[$\ldots$]
		[$\ldots$]] 
		[
		{$(\overline{(2,4)_2,2,3,4})$}
		[$\ldots$]
		[$\ldots$]]]]]
	\end{forest}
	\caption{The Markov-Hurwitz tree $\mathcal{MH}$}
\end{figure}

 We define the 
\textit{conjunction operation} of two periods as
\begin{equation}\label{conj}
	(\overline{a_1,\ldots,a_r})\odot(\overline{b_1,\ldots,b_s})=(\overline{a_1,\ldots,a_r,b_1,\ldots,b_s}).
\end{equation}
Each Markov irrationality not on  the most left branch is 
 the result of the conjunction operation of two predecessors: its immediate predecessor on the same branch and the immediate predecessor 
of the tip of the branch. We call \textbf{right} predecessor the one on the right of $w$ on the tree, and \textbf{left} predecessor the one on the left. With this terminology, $w$ is the result of the cojunction operation of the right predecessor with the left predecessor. For example, 
$$(\overline{2,4,(2,3,4)_2}) = (\overline{2,4,2,3,4})\cdot (\overline{2,3,4});
$$ 
the right and left predecessors are $(\overline{2,4,2,3,4})$ and $(\overline{2,3,4})$ respectively. On the most left  branch, the Markov irrationality at level $n$ is $(\overline{2,3_n,4})$. We call $(\overline{3})$ the \textbf{left} predecessor and $(\overline{2,3_{n-1},4})$ the \textbf{right} predecessor.

\subsection{The Farey tree}

  There is a natural parametrisation of  Markov numbers by  Farey fractions which goes back to Frobenius \cite{fro}. 
  The Farey tree is constructed by following a very similar procedure as described in section 3.1 for  the Markov-Hurwitz tree.
   We start with $\frac{0}{1}$, $\frac{1}{2}$ and $\frac{1}{3}$, and construct  
each following fraction from the two predecessors $\frac{a}{b}, \frac{c}{d}$ that are in the same position as in the Markov-Hurwitz tree ($\frac{a}{b}$ is the right predecessor and $\frac{c}{d}$ is the left predecessor) by taking the Farey median 
$$
\dfrac{a}{b}*\dfrac{c}{d}=\dfrac{a+c}{b+d}.
$$
Hence  each vertex of the tree is a Farey fraction $p/q$ that is in correspondence with  a 
 Markov quadratic $w=w(p/q)$ and a Markov number $c = m(p/q)$. We denote the Farey tree by $\F$.

\begin{figure}[h!]
	\centering
\begin{forest}
[
{$\dfrac{0}{1}$\qquad$\dfrac{1}{2}$}
[
{$\dfrac{1}{3}$} 
[
{$\dfrac{1}{4}$} 
[
{$\dfrac{1}{5}$}
[]
[]]
[
{$\dfrac{2}{7}$}
[]
[]
]] 
[
{$\dfrac{2}{5}$} 
[
{$\dfrac{3}{8}$}
[]
[]] 
[
{$\dfrac{3}{7}$}
[]
[]]]]]
\end{forest}
\caption{The Farey tree $\F$}
\end{figure}
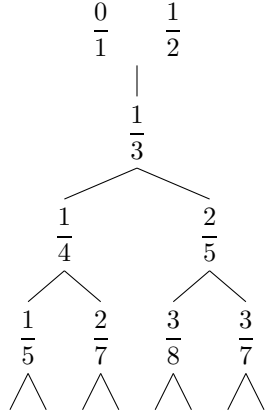

 By the construction of the Markov-Hurwitz and the Farey trees, the  denominators in $\F$  correspond to the  lengths of the periods of the continued fractions  that occupy the same corresponding position in  $\mathcal{MH}$.

The operation giving the denominators in the Farey tree corresponds to  the well known Euclidean algorithm which, starting with  $(0,1,1)$, gives from a triple $(s,t,u)$
two new triples
 $(s,u,s+u)$ and $(t,u,t+u)$, obtaining thus all solutions 
 to the equation 
\begin{equation}\label{Euclideq}
s+t=u, \quad 0\leq s\leq t\leq u,\quad (s,t)=1.
\end{equation}
The denominators in the Farey tree are the maximum terms of the solutions $(s,t,u)$.
 
\section{Asymptotic formula for $\log\varepsilon$ in terms of Farey denominators}


We order the denominators of the Farey fractions in $\F$  as real increasing numbers $(q_n)_{n\geq1}$. If two Farey fractions $\frac{p}{q},\frac{p'}{q'}$ satisfy $q=q'$, then we write $q=q_n$, $q'=q_m$ with $n<m$ if $m(\frac{p}{q})<m(\frac{p'}{q'})$. Hence we obtain a sequence:    
\begin{equation}\label{qn}
	q_1= 1,\quad q_2=2,\quad q_3=3,\quad q_4=4, \quad q_5=5,\quad q_6=5,\ldots
\end{equation}
that is in bijective correspondence with the sequence of Markov numbers \eqref{markovseq}.
\\

The number of denominators $q_n$ less than $x$ (with multiplicity) is given by
$$
 \#\left\{q_n\leq x\right\}_{n\geq 1} = 
 \# \left\{(s,t,u)\ solution\ to\ \eqref{Euclideq}:\, u\leq x\right\}=1+\dfrac{1}{2}\sum_{u\leq x} \varphi(u),
$$
where  $\varphi$ is the  Euler function. From the
 asymptotic formula 
$$1+\dfrac{1}{2}\sum_{u\leq x} \varphi(u) \sim\dfrac{3}{2\pi^2} x^2,$$ 
 it follows that
\begin{equation}\label{asimqn}
	q_n\sim \pi\sqrt{\frac{2}{3}n}.
\end{equation}

Let $w=w(p/q)$ be a Markov irrationality and $p/q$ be  its corresponding Farey fraction.
Set $q=q_n$, where $q_n$ is a Farey denominator ordered as in \eqref{qn}, so that $m(p/q)=c_n$.
The discriminant of $w$ is $9c_n^2-4$, so Pell's equation is
$$
t^2-(9c_n^2-4)u^2=4,
$$
with $(3c_n,1)$ being the smallest positive solution. Thus $\varepsilon=\frac{3c_n+\sqrt{9c_n^2-4}}{2}$.

By a result of Zagier \cite{ZMar}, we have that
\begin{equation}\label{asimcn}
	c_n\sim e^{\sqrt{\frac{n}{C}}},
\end{equation}
where $C\approx 0.18071704711507$.
Using \eqref{asimqn},
\begin{equation}\label{asycn}
	c_n\sim e^{\frac{q_n}{\pi}\sqrt{\frac{3}{2C}}}.
\end{equation}
Hence
\begin{equation}\label{loge}
	\log\varepsilon\sim \dfrac{\sqrt{3}}{\sqrt{2C}\pi} q_n + \log\frac{3}{2}. 
\end{equation}

\section{Simple forms and continued fractions}
 
Any indefinite binary quadratic form $[a,b,c]$ is $\SL$-equivalent to one satisfying
 \begin{equation} 	\label{simple}
 a>0>c.
 \end{equation}
There is a finite number of forms satisfying \eqref{simple} with fixed discriminant. Such forms are commonly called `simple' after Zagier \cite{Zbook} and  play an important role in the theory of rational periods and  period functions of modular forms (see \cite{KZ84}, \cite{Z99}, \cite{Ch}, \cite{PB-JNT}).  
 
 It is shown in \cite{ChZ} that all simple  forms  $\SL$ equivalent to a form $[a,b,c]$ are obtained by applying iteratively the following  continued fraction algorithm.  Let $w$ be the period  in the `--' continued fraction expansion of the root $\frac{-b+\sqrt{b^2-4ac}}{2a}$; we define
 \begin{equation}\label{wi}
w^{(1)}=w-1,\qquad
 w^{(i+1)}=\left\{\begin{array}{ll}
w^{(i)}-1=\tiny{\begin{pmatrix} 1 &-1\\0 &1\end{pmatrix}} (w^{(i)}) &\mbox{if $w^{(i)}\geq 1 $},\\
\dfrac{w^{(i)}}{1-w^{(i)}}=\tiny{\begin{pmatrix} 1 &0\\-1 &1\end{pmatrix}} (w^{(i)}) &\mbox{otherwise}.
\end{array}\right.
\end{equation}
 This algorithm is cyclic because the continued fraction of $w$ is purely periodic; we denote by $\l=\l_w$ the length of the cycle, 
 so that $w^{(\l+1)}=w^{(1)}$.
 The cycle $w^{(1)},\ldots,w^{(\l)}$ corresponds to the cycle of simple forms  in the $\SL$-equivalence class of the quadratic form  $[a,b,c]$. 
 
In \cite{BI1} the value $j(w)$ is written in terms of the cycle \eqref{wi} when $w$ is a Markov irrationality (see \cite[Lemma 4.1]{BI1}) and in \cite{BI2} for arbitrary quadratic irrationalities. More concretely, in terms of simple forms, we have:
 \begin{lemma}\label{lemma2} For any real quadratic irrationality $w$,
\begin{equation}\label{**}
j(w)= \dfrac{1}{2\log\varepsilon} \int_{e^{\pi i/3}}^{e^{2\pi i/3}} j(z) \sum_{\substack{\left[a,b,c\right]\ \mathrm{simple}\\ \left[a,b,c\right]\in\A}} \dfrac{\sqrt{D}}{az^2+bz+c} dz,
\end{equation} 
where $\A$ is the $\SL$-equivalence class of $w$ and $D$ its discriminant. 
\end{lemma}
If $w$ has a purely periodic `--' continued fraction expansion, we can write \eqref{**}  as:
\begin{equation}	\label{fkernel}
	j(w)=  \dfrac{1}{2\log\varepsilon} \int_{\pi/3}^{2\pi /3} j(e^{i\theta}) ie^{i\theta} \sum_{i=1}^{\l_w} \Big(\dfrac{1}{e^{i\theta}-w^{(i)}}-\dfrac{1}{e^{i\theta}-\tilde w^{(i)}}\Big) d\theta,
\end{equation}
where  $w^{(i)}$ are defined in (\ref{wi}) and $\tilde{w}^{(i)}$ are  the Galois conjugates of $w^{(i)}$. 
Each term $w^{(i)}$ is of the form 
$$
w^{(i)}=(a_0,\overline{a_1,\ldots,a_n}) = a_0-\dfrac{1}{(\overline{a_1,\ldots,a_n})}\qquad (1\leq a_0\leq a_n-1).
$$
It is a well known fact (see for example \cite{Zbook}) that the conjugate of $1/(\overline{a_1,\ldots,a_n})$ is $(\overline{a_n,\ldots,a_1})$ and hence 
$$
\tilde w^{(i)}= -(a_n-a_0,\overline{a_{n-1},a_{n-2},\ldots,a_1,a_n}).
$$
If $w$ is a Markov irrationality, the structure of the Markov-Hurwitz tree implies that:
\begin{align}	 
	\dfrac{3}{8}=(1,2,3,2) \leq &w^{(i)}\leq (3,2,4,2)=\dfrac{29}{12},\label{wibounds}\\
	-\dfrac{21}{8}=-(3,3,3)\leq &\tilde w^{(i)} \leq  -(1,2,3)=-\dfrac{2}{5}.\label{wibounds2} 
\end{align}

\section{Stratregy for the proof of Theorem \ref{th1}}

Let $w$ be a Markov irrationality on the $n$-th level of the Markov-Hurwitz tree, the first level corresponding to the Markov irrationality $(\overline{2,3,4})$.
 Let $v$ be the `immediate' predecessor of $w$, which is on the same branch as $w$. Then we have one of the following two configurations:
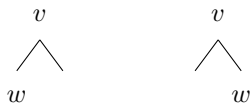
\begin{center}
	\begin{forest}
		[{$v$}
		[{$w$}]
		[{$ $}]]
		\end{forest} \qquad\qquad
	\begin{forest}
		[{$v$}
		[{$ $}]
		[{$w$}]]
		\end{forest}
	\end{center}
Let  $u=(\overline{a_1,\ldots,a_s})$ be the other predecessor of $w$ as explained in section 3.1. In the second configuration, $u$ is the right predecessor of $w$ and $v$. Then let $(\overline{b})=(\overline{b_1,\ldots,b_t})$ be the left predecessor of $u$. We have  $v=(\overline{u_r,b})$ for some  $r\geq 1$ and $w=(\overline{u_{r+1},b})$. 
In the first configuration, $u$ is the left predecessor of $w$ and $v$. 
If $w$ is not on the most left branch of the tree, then we set $(\overline{b})=(\overline{b_1,\ldots,b_t})$ be the right predecessor of $u$, so that $v=(\overline{b,u_r})$ for some $r\geq 1$ and $w=(\overline{b,u_{r+1}})$. The cycles are the same if we consider $(\overline{u_r,b})$ and $(\overline{u_{r+1},b})$, and we will refer to these two situations  as \textbf{`case 1'}.
Note that in this case, $rs>n$.

The situation where  $w$ is  on the most left branch of the tree will be  \textbf{`case 2'}. In this case, $u=(\overline{3})$, $v=(\overline{2,3_r,4})$ for some $r\geq 1$ and $w=(\overline{2,3_{r+1},4})$.
Note that in this case, $r=n-1$ and
$$
\ell_u=2,\quad\ell_v=2r+4,\quad \ell_w=\ell_u+\ell_v=2(r+1)+4.
$$

Let 
$$
J(w):= (2\log\varepsilon) j(w).
$$
In the next section we will write $J(w)$ as
\begin{equation}\label{Jstra}
J(w)=J(u)+J(v) +\delta_w,
\end{equation}
where $\delta_w$ is the error term that we will explicitly bound in terms of $n$. We see the identity \eqref{Jstra} as a `local' formula for $J(w)$ on the Markov-Hurwitz tree.
We will deduce a recursive formula for $J(w)$ and $q_n$ in section 8 (see \eqref{recJ2}, \eqref{recJ3}, \eqref{recq2}, \eqref{recq3}) that we see as a `global' formula and that will give the global interlacing property in section 9.

In order to obtain the local formula \eqref{Jstra}  we will compare the cycle of $w$ in the algorithm \eqref{wi} with the cycles of $u$  and $v$ and their $J$ values in section 7.1. In sections 7.2, 7.3, 7.4 we bound the error term $\delta_w$.

\section{Local formula for $(2\log\varepsilon) j(w)$}

\subsection{Comparing cycles of $w,u,v$.} We will compare the cycle of $w$ in the algorithm \eqref{wi} with the cycles of $u$  and $v$. We define a sum $S_u$ and a sum $S_v$ that contain all the terms in the cycle of $u$ and $v$ respectively. The definitions of $S_u$ and $S_v$ are different for case 1 and case 2. In case 1, 
the sum $S_u$ below compares the first  terms $w^{(1)},\ldots,w^{(\ell_u)}$ and the conjugates $\tilde w^{(r\ell u+1)},\ldots,\tilde w^{((r+1)\ell_u)}$ with the cycle of $u$ and  conjugates, while the sum $S_v$ compares the remaining terms  with the cycle of $v$ and conjugates. 
\\

\noindent
\textbf{Case 1:}
\begin{align*}
	&S_u(\theta)=\sum_{i=1}^{\ell_u} \dfrac{1}{e^{i\theta}-w^{(i)}} - \dfrac{1}{e^{i\theta}-u^{(i)}}+
	\sum_{i=1}^{\ell_u} \dfrac{1}{e^{i\theta}-\tilde u^{(i)}} -  \dfrac{1}{e^{i\theta}-\tilde w^{(r\ell_u+i)}}; \\
	&S_v(\theta)=\sum_{i=1}^{\ell_v} \dfrac{1}{e^{i\theta}-w^{(i+\ell_u)}}-\dfrac{1}{e^{i\theta}-v^{(i)}}
	+\sum_{i=1}^{r\ell_u}   \dfrac{1}{e^{i\theta}-\tilde v^{(i)}} - \dfrac{1}{e^{i\theta}-\tilde w^{(i)}}\\
	& + \sum_{i=r\ell_u+1}^{\ell_v}   \dfrac{1}{e^{i\theta}-\tilde v^{(i)}} - \dfrac{1}{e^{i\theta}-\tilde w^{(\ell_u+i)}}.
\end{align*}

\textbf{Case 2:}
\begin{align*}
	&S_u(\theta)=\sum_{i=1}^{2} \dfrac{1}{e^{i\theta}-w^{(1+i)}} - \dfrac{1}{e^{i\theta}-u^{(i)}}+\sum_{i=1}^2  
	\dfrac{1}{e^{i\theta}-\tilde u^{(i)}}
	-\dfrac{1}{e^{i\theta}-\tilde w^{(1+2r+i)}};\\
	&S_v(\theta)=\dfrac{1}{e^{i\theta}-w^{(1)}} -\dfrac{1}{e^{i\theta}-v^{(1)}}+  \sum_{i=2}^{4+2r} \dfrac{1}{e^{i\theta}-w^{(2+i)}}-\dfrac{1}{e^{i\theta}-v^{(i)}}\\
	&+ \sum_{i=1}^{1+2r} \dfrac{1}{e^{i\theta}-\tilde v^{(i)}} - \dfrac{1}{e^{i\theta}-\tilde w^{(i)}}  + \sum_{i=2+2r}^{4+2r}  
	\dfrac{1}{e^{i\theta}-\tilde v^{(i)}}
	- \dfrac{1}{e^{i\theta}-\tilde w^{(2+i)}}.
\end{align*}

Using \eqref{fkernel}, we have that
\begin{equation}\label{intwithsums}
	J(w) - J(u)- J(v)=  \int_{\pi/3}^{2\pi/3} j(e^{i\theta}) i e^{i\theta} (S_u(\theta) + S_v(\theta)) d\theta.
\end{equation}
Our next goal is to bound the integral in \eqref{intwithsums}; we will give bounds for the real and imaginary parts.
For $\alpha\in\left\{u,v\right\}$, let 
\begin{align*}
	\varepsilon_\alpha(\theta)&=\cos\theta\Im(S_\alpha(\theta)) + \sin\theta\Re(S_\alpha(\theta));\\
	\varepsilon'_\alpha(\theta)&=\cos\theta\Re(S_\alpha(\theta)) - \sin\theta\Im(S_\alpha(\theta)).
\end{align*}
The real and imaginary parts of the integral in \eqref{intwithsums} are respectively
\begin{equation}	\label{intreal}
\Re(\eqref{intwithsums})=-\int_{\pi/3}^{2\pi/3} j(e^{i\theta}) (\varepsilon_u(\theta)+\varepsilon_v(\theta)) d\theta
\end{equation}
and
\begin{equation}	\label{intimag}
\Im(\eqref{intwithsums})=\int_{\pi/3}^{2\pi/3} j(e^{i\theta}) (\varepsilon'_u(\theta)+\varepsilon'_v(\theta)) d\theta.
\end{equation}
Our goal is to bound \eqref{intreal} and \eqref{intimag}.
It will be useful to define the following real functions with three variables $(x,y,\theta)\in\mathbb{R}^2\times[\frac{\pi}{3},\frac{2\pi}{3}]$:
$$	
g(x,y,\theta)=	\dfrac{1}{x-y}\Big(\cos\theta\Im\Big(\dfrac{1}{e^{i\theta}-x} - \dfrac{1}{e^{i\theta}-y}\Big) + \sin\theta\Re\Big(	\dfrac{1}{e^{i\theta}-x} - \dfrac{1}{e^{i\theta}-y}\Big)\Big),
$$
$$	
g'(x,y,\theta)= \dfrac{1}{x-y}	\Big(\cos\theta\Re\Big(\dfrac{1}{e^{i\theta}-x} - \dfrac{1}{e^{i\theta}-y}\Big) - \sin\theta\Im\Big(	\dfrac{1}{e^{i\theta}-x} - \dfrac{1}{e^{i\theta}-y}\Big)\Big).
$$
Note that
\begin{equation}\label{g}
	\dfrac{1}{e^{i\theta}-x} - \dfrac{1}{e^{i\theta}-y} = \dfrac{(x-y)(\Re((e^{i\theta}-x)(e^{i\theta}-y))-i\Im((e^{i\theta}-x)(e^{i\theta}-y))}{|(e^{i\theta}-x)(e^{i\theta}-y)|^2}.
\end{equation}
We have that
\begin{align*}
	&\Re((e^{i\theta}-x)(e^{i\theta}-y))= \cos^2\theta - \sin^2\theta - \cos\theta(x+y)+xy,\\
	&\Im((e^{i\theta}-x)(e^{i\theta}-y))=2\sin\theta\cos\theta-\sin\theta(x+y),
\end{align*}
so
\begin{align}
	g(x,y,\theta) &= -\dfrac{\sin\theta(1-xy)}{((\cos\theta-x)^2+\sin^2\theta)((\cos\theta-y)^2+\sin^2\theta)}\label{real}\\
	g'(x,y,\theta)&=\dfrac{ -x -y +\cos\theta(1+xy)}{((\cos\theta-x)^2+\sin^2\theta)((\cos\theta-y)^2+\sin^2\theta)}. \label{imag}
\end{align}
For  $ x,y$ satisfying \eqref{wibounds}, one can easily compute:
\begin{align}\label{num1}
	-1.26964  \leq	g(x,y,\theta) \leq 0.354112 ,\quad
	-1.10636\leq	g'(x,y,\theta) \leq -0.07222  .
\end{align}
For $x,y$ satisfying \eqref{wibounds2},
\begin{align}\label{num3}
	-1.25946\leq	g(x,y,\theta) \leq 0.354112, \quad
	0.04705	\leq g'(x,y,\theta) \leq  1.10636.
\end{align}
The following lemma (see \cite[Lemma 2.1]{BI1} for a proof) will also be useful:
\begin{lemma}\label{lemacoincide}
	If the `--' continued fraction expansions of two Markov quadratics $u$ and $v$ coincide in the first $r$ partial quotients, then 
	$$
	|u-v|\leq 10\Big(\dfrac{2}{1+\sqrt{5}}\Big)^{2(r-1)}.
	$$
\end{lemma}

Let
$$
b(x)=10\Big(\dfrac{2}{1+\sqrt{5}}\Big)^{2(x-1)}
$$
be the bound of Lemma \ref{lemacoincide}. Note that
 \begin{equation}\label{sumb}
 	\sum_{k=k_0}^\infty b(k)  \leq 10\Big(\dfrac{2}{1+\sqrt{5}}\Big)^{2k_0-3}.
 \end{equation}

\subsection{Bounds for $\varepsilon_u,\varepsilon'_u$.}

\noindent
\textbf{Case 1.}
We write  $\varepsilon_u$ using the function $g(x,y,\theta)$:
$$
\varepsilon_u(\theta)= \sum_{i=1}^{\ell_u} (w^{(i)}-u^{(i)}) g(w^{(i)},u^{(i)},\theta) + \sum_{i=1}^{\ell_u} (\tilde u^{(i)}-\tilde w^{(r\ell_u+i)}) g(\tilde u^{(i)},\tilde w^{(r\ell_u+i)}, \theta).
$$
We have a similar expression for $ \varepsilon'_u(\theta)$ replacing $g$ by $g'$.
We apply Lemma \ref{lemacoincide} to bound $|w^{(i)}-u^{(i)}|$ and $|\tilde u^{(i)}- \tilde w^{(r\ell_u+i)}|$ and use  \eqref{num1}, \eqref{num3} to bound the rest.

In the first sum of $\varepsilon_u(\theta)$, 
the first $a_1-1$ terms $w^{(i)}$ and $u^{(i)}$ share the same first $rs+s$ partial quotients,  the next $a_2-1$ terms share the first $rs+s-1$ partial quotients, etc. and the last $a_s-1$ terms coincide in at least the first $rs+1$ partial quotients.

In the second sum $\varepsilon_u(\theta)$, we have a similar reversed situation, where the first $a_1-1$ terms $\tilde w^{(r\ell_u+i)}$ and $\tilde u^{(i)}$ share the same first $rs+1$ partial quotients, the next  $a_2-1$ terms share the same $rs+2$ first partial quotients, etc. and the last $a_s-1$ terms share the same $rs+s$ partial quotients. 


Using  Lemma \ref{lemacoincide},  \eqref{num1}, \eqref{num3} and that  $a_i-1\leq 3$ for $i=1,\ldots,s$, we have
\begin{align*}\label{intS1}
	&|\varepsilon_u(\theta)| \leq   3(1.26964 +1.25946)\,  \sum_{k=1}^{s} b(rs+k)  \\ 
	&|\varepsilon'_u(\theta)|\leq   1.10636 \cdot 6 \sum_{k=1}^{s} b(rs+k).
\end{align*}

\noindent
\textbf{Case 2.} We proceed analogously to case 1.
The two terms in each  sum of $S_u$ share the same  first $r+1$ partial quotients. By Lemma \ref{lemacoincide} we have, for $i=1,2$,
\begin{equation}	\label{case2}
|w^{(1+i)}-u^{(i)}|, |\tilde w^{(i+1+2r)}-\tilde u^{(i)}|\leq b(r+1).
\end{equation}
Using    \eqref{num1}, \eqref{num3}  and \eqref{case2}, we have that
\begin{align*}
	&|\varepsilon_u(\theta)| \leq   2\cdot 2.5291\, b(r+1)\\ 
	&|\varepsilon'_u(\theta)|\leq   4\cdot1.10636\, b(r+1).
\end{align*}

\subsection{Bounds for $\varepsilon_v,\varepsilon_v'$.} 
We bound $\varepsilon_v,\varepsilon_v'$ in a similar way to $\varepsilon_u,\varepsilon_u'$.

\noindent
\textbf{Case 1.}
In the first sum of $S_v$,
the first $a_1-1$ terms $w^{(\ell_u+i)}$ and $v^{(i)}$ coincide in the first  $2rs+t$ partial quotients, the next $a_2-1$ terms share the same first $2rs+t-1$ partial quotients, etc. and the last $b_t-1$ terms coincide in the first  $rs+1$ partial quotients.

In the second sum of $S_v$,  the first $a_1-1$ terms $\tilde w^{(i)}$ and $\tilde v^{(i)}$  share the same first $rs+t+1$ partial quotients, the next $a_2-1$ terms coincide in the first $rs+t+2$ partial quotients, etc. and the last $a_s-1$ terms coincide in the first $2rs+t$ partial quotients. 

In the third sum of $S_v$, the first $b_1-1$ terms $\tilde v^{(i)}$ and $\tilde w^{(\ell_u+i)}$ coincide in the first $rs+1$ partial quotients, the next $b_2-1$ terms coincide in the first $rs+2$ partial quotients, etc. and the last $b_t-1$ partial quotients coincide in the first $rs+t$ partial quotients.

Using again Lemma \ref{lemacoincide},  \eqref{num1}, \eqref{num3} and that  $a_i-1\leq 3$, $b_j-1\leq 3$ for $i\in\left\{1,\ldots,s\right\}, j\in\left\{1,\ldots,t\right\}$, we have
\begin{align*}
	&|\varepsilon_v(\theta)|\leq  3\cdot 2.5291\,  \sum_{k=1}^{rs+t} b(rs+k)\\
	&	|\varepsilon'_v(\theta)| \leq  1.10636 \cdot 6 \sum_{k=1}^{rs+t} b(rs+k).
\end{align*}

\noindent
\textbf{Case 2.}
The terms $w_1$ and $v_1$ share  the first same $r+1$ partial quotients.
In the first sum of $S_v$, the first two terms coincide in $2r+2$ partial quotients, each next two terms coincide in $2r+1, \ldots, r+3$ partial quotients, and the last three terms coincide in $2+r$ partial quotients.    
In the second sum, the first terms share the same $1+r$ first partial quotients, and each next block of two terms $\tilde w_i, \tilde v_i$  coincide in $2+r$, $3+r$, etc.. $2+2r$  first partial quotients. In the third sum, the terms coincide in $r+1$ partial quotients.  

Using  Lemma \ref{lemacoincide}, \eqref{num1}, \eqref{num3}, we have that
\begin{align*}
	&|\varepsilon_v(\theta)| \leq   1.26964 (b(r+1) + b(r+2)) +2.5291 \cdot 2\sum_{k=2}^{r+2} b(r+k) + 1.25946\cdot 4b(r+1)\\ 
	&|\varepsilon'_v(\theta)|\leq   1.10636\Big(5b(r+1) + b(r+2) +4\sum_{k=2}^{r+2} b(r+k) \Big).
\end{align*}

\subsection{Conclusion.}  Recall that we denote by $n$ the level on the tree where the Markov irrationality $w$ lies.
In case 1,  using that $rs>n$, we have
\begin{align}
|\varepsilon_u(\theta)| + |\varepsilon_v(\theta)|
&\leq 151.7460\Big(\dfrac{2}{1+\sqrt{5}}\Big)^{2n+1},
 \label{concl1}\\
|\varepsilon'_u(\theta)| + |\varepsilon'_v(\theta)| 
&\leq  132.7632\Big(\dfrac{2}{1+\sqrt{5}}\Big)^{2n+1}.
 \label{concl2}
\end{align}
In case 2, using that $r=n-1$, we have 
\begin{align}
|\varepsilon_u(\theta)| + |\varepsilon_v(\theta)| &\leq 
113.6564 \Big(\dfrac{2}{1+\sqrt{5}}\Big)^{2(n-1)}  +63.278 \Big(\dfrac{2}{1+\sqrt{5}}\Big)^{2n-1}, 
 \label{concl3}\\
|\varepsilon'_u(\theta)| + |\varepsilon'_v(\theta)| &\leq 
99.5724\Big(\dfrac{2}{1+\sqrt{5}}\Big)^{2(n-1)}+ 55.318 \Big(\dfrac{2}{1+\sqrt{5}}\Big)^{2n-1}.  \label{concl4}
\end{align}
Using \eqref{concl1}, \eqref{concl2}, \eqref{concl3}, \eqref{concl4} and that 
\begin{equation*}	\label{average}
\int_{\pi/3}^{2\pi/3} j(e^{i\theta}) d\theta=753.982,
\end{equation*}
we have
\begin{equation}	\label{ind}
J(w)=J(u) +J(v) +\delta_w,
\end{equation}
with 
\begin{align}
&|\Re(\delta_w)|\leq  115181.57371 \Big(\dfrac{2}{1+\sqrt{5}}\Big)^{2(n-1)} \label{bound delta1}\\
&|\Im(\delta_w)|\leq  100853.23866\Big(\dfrac{2}{1+\sqrt{5}}\Big)^{2(n-1)}. \label{bound delta2}
\end{align}

\section{Global formula for $(2\log\varepsilon) j(w)$}

Consider a path on the Markov-Hurwitz tree given by Markov irrationalities $w_0$, $w_1$, $w_2$, $\ldots, w_{n}$, $n\geq2$.  The quadratic $w_n$ is on the $n$-th level of the tree,  the first level corresponding to  $w_1=(\overline{2,3,4})$. 
If $w_n$ is on the left half of the tree, $w_0=(\overline{3})$, otherwise $w_0=(\overline{2,4})$.
 Let $r_1,\ldots,r_m$ be the successive levels where the path turns, starting from $r_1=1$, and let $r_{m+1}=n$. For convenience we set $r_{-1}=1$.
  We write $q_i$ for the denominator of the Farey fraction associated to $w_i$. If $w_n$ is on the left half of the tree, then $q_0=1$, otherwise $q_0=2$. We write $\delta_i=\delta_{w_i}$ for the error in \eqref{ind} associated to $w_i$.

It follows from the relation \eqref{ind} applied to $J(w_n)$, that
\begin{equation}	\label{Jwn1}
J(w_n)=J(w_{n-1}) + J(w_{r_m-1}) + \delta_n.
\end{equation}
By the construction of the Farey tree, we also have
\begin{equation}	\label{qn1}
q_n=q_{n-1}+q_{r_m-1}.
\end{equation}

If we apply \eqref{ind} recursively to  $J(w_{n}),\dots,J(w_{r_m})$ and $q_n,\ldots,q_{r_m}$ we have, if $m\geq 2$,
\begin{equation}	\label{recJn}
	J(w_n)= (n-r_m+1) J(w_{r_m-1}) + J(w_{r_{m-1}-1})+ \sum_{i=r_m}^n\delta_i,
\end{equation}
\begin{equation}	\label{recqn}
	q_n=(n-r_m+1) q_{r_m-1}+  q_{r_{m-1}-1}.
\end{equation} 
If we apply again \eqref{ind} recurrently to $J(w_{r_k-1}),\dots,J(w_{r_{k-1}})$ and $q_{r_k-1}, \ldots,q_{r_{k-1}}$ we obtain,  for $3\leq k\leq m$:
\begin{equation}	\label{recJ}
J(w_{r_k-1})= (r_k-r_{k-1}) J(w_{r_{k-1}-1}) + J(w_{r_{k-2}-1})+ \sum_{i=r_{k-1}}^{r_k-1}\delta_i,
\end{equation}
\begin{equation}	\label{recq}
q_k=(r_k-r_{k-1}) q_{r_{k-1}-1}+  q_{r_{k-2}-1}.
\end{equation}

Let $k_0$ be a positive constant `large enough' that will indicate a level `down enough' on the Markov-Hurwitz tree. We will choose a specific value for $k_0$ later. 
We denote by $k$ and $s$ the positive integers such that
$$
k_0=r_{k-1} + s,\qquad  1\leq s\leq r_{k}-r_{k-1}.
$$
The recursive formulas \eqref{recJ},\eqref{recJn} give, if $m\geq 2$,
\begin{equation}	\label{recJ2}
\begin{pmatrix}
	J(w_n)\\ J(w_{r_m-1})
\end{pmatrix}
=
\begin{pmatrix} n-r_m+1 &1\\
	1 &0\end{pmatrix}
\Big(\prod_{i={k+1}}^m \begin{pmatrix} r_i-r_{i-1} &1\\
	1 &0\end{pmatrix}\Big) \begin{pmatrix} s &1\\1 &0\end{pmatrix}
\begin{pmatrix}
	J(w_{r_{k-1}-1})\\ J(w_{r_{k-2}-1})
\end{pmatrix}
+\begin{pmatrix}
	\delta \\ \delta'
\end{pmatrix}
\end{equation}
where 
\begin{equation}	\label{deltafin}
\delta= \sum_{i=r_m}^n \delta_i+ \Big( \sum_{j=k+1}^{m} \lambda_j \sum_{i=r_{j-1}}^{r_{j}-1} \delta_i \Big) +  \lambda_k \sum_{i=k_0}^{r_{k}-1} \delta_i
\end{equation}
and $\lambda_j$ is the coefficient of $J(w_{r_{j-1}})$ in the recursive formula \eqref{recJ2} for $k\leq j\leq m$.

If $m=1$, then
\begin{equation}\label{recJ3}
J(w_n)=(n-k_0) J(w_{0}) + J(w_{k_0}) + \sum_{i=k_0-1}^n \delta_i.
\end{equation}
We obtain similar formulas for $q_n$ from \eqref{recq},\eqref{recqn}: if $m\geq 2$,
\begin{equation}	\label{recq2}
\begin{pmatrix}
	q_n\\ q_{r_m-1}
\end{pmatrix}
=
\begin{pmatrix} n-r_m+1 &1\\
	1 &0\end{pmatrix}
\Big(\prod_{i={k+1}}^m \begin{pmatrix} r_i-r_{i-1} &1\\
	1 &0\end{pmatrix}\Big)\begin{pmatrix} s &1\\1 &0\end{pmatrix}
\begin{pmatrix}
	q_{r_{k-1}-1}\\ q_{r_{k-2}-1}
\end{pmatrix}
\end{equation}
and if $m=1$,
\begin{equation}	\label{recq3}
q_n=(n-k_0) q_{0} + q_{k_0}.
\end{equation}
The coefficients of  $q_{r_{j-1}}$ in \eqref{recq2} are the same as the coefficients $\lambda_j$ for $J(w_{r_{j-1}})$ in \eqref{recJ2}, so we have, 
for all $j\geq 2$, 
\begin{equation}	\label{qnbig}
q_n\geq \lambda_j q_{r_{j-1}}
\end{equation}
and
\begin{equation}\label{Jn}
	\dfrac{J(w_{n})}{q_{n}}=
	\dfrac{\alpha J(w_{r_{k-2}-1}) + \beta J(w_{r_{k-1}-1})}{\alpha q_{r_{k-2}-1}  +  \beta q_{r_{k-1}-1}} + \dfrac{\delta}{q_{n}},
\end{equation}
where $\alpha,\beta$ are some positive constants given by \eqref{recJ2}, \eqref{recq2} and $\delta$ is given by \eqref{deltafin}.

\section{Asymptotic interlacing property}

Suppose $n>k_0$.
Below we bound  $\delta/q_n$. 
We consider first the case when $m\geq 2$.
Using  \eqref{bound delta1}, \eqref{qnbig} and $q_{r_{k-1}}\geq 4$, we have that 
\begin{align}\label{erdeltak}
	\dfrac{\lambda_k}{q_n} \sum_{i=k_0}^{r_k-1} |\Re(\delta_i)| &\leq 
		\dfrac{115181.57371}{4} \sum_{i=k_0}^{r_k-1} \Big(\dfrac{2}{1+\sqrt{5}}\Big)^{2(i-1)}\nonumber\\
		&\leq  \dfrac{115181.57371}{4} \Big(\dfrac{2}{1+\sqrt{5}}\Big)^{2k_0-3}.
	\end{align}
Similarly, for $j\geq k+1$, using \eqref{qnbig} and $q_{r_j-1}\geq r_{j-1}\geq k_0$, 
\begin{align}\label{erdeltaj}
\dfrac{1}{q_n}  \sum_{j=k+1}^{m} \lambda_j  	 \sum_{i=r_{j-1}}^{r_j-1} |\Re(\delta_i)|
&\leq 
\sum_{j=k+1}^m \dfrac{1}{q_{r_{j-1}}}  \sum_{i=r_{j-1}}^{r_j-1} |\Re(\delta_i)| \nonumber\\
&\leq \dfrac{115181.57371}{k_0} \Big(\dfrac{2}{1+\sqrt{5}}\Big)^{2k_0-1}.
\end{align}
By \eqref{bound delta1} and $q_n> k_0$, we also have that
\begin{equation}\label{erdeltan}
\dfrac{|\Re(\delta_n)|}{q_n}\leq  \dfrac{115181.57371 }{k_0+1}  \Big(\dfrac{2}{1+\sqrt{5}}\Big)^{2k_0}. 
\end{equation}
If $m=1$, then $q_n\geq n+2 \geq k_0+3$, so
\begin{equation}\label{deltam1}
\dfrac{|\Re(\delta)|}{q_n} = \dfrac{1}{q_n} \sum_{i=k_0-1}^n  |\Re(\delta_i)| \leq \dfrac{115181.57371}{k_0+3}    \Big(\dfrac{2}{1+\sqrt{5}}\Big)^{2k_0-5}.
\end{equation}

Similar bounds for the imaginary parts can be found by using \eqref{bound delta2}. Thus we have, letting $k_0\rightarrow\infty$ as $n\rightarrow\infty$, 
\begin{equation}	\label{deltainf}
\frac{\delta}{q_n}
\rightarrow 0\qquad as \ n\rightarrow\infty.
\end{equation}

By \eqref{Jwn1} and \eqref{qn1} we have, if $\frac{\Re(J(w_{r_m-1}))}{q_{r_m-1}}\leq  \frac{\Re(J(w_{n-1}))}{q_{n-1}}$, 
$$
\dfrac{\Re(J(w_{r_m-1}))}{q_{r_m-1}}\leq \dfrac{\Re(J(w_n))}{q_n}\leq \dfrac{\Re(J(w_{n-1}))}{q_{n-1}}
\qquad as \ n\rightarrow\infty,
$$
and otherwise the inequalities are reversed. The same inequalities for the imaginary parts hold.
By  \eqref{asimqn}, we conclude that
$j(w_n)$ lies between $j(w_{n-1})$ and $j(w_{r_{m}-1})$, and hence Theorem \ref{th1} is proved.

\section{Asymptotic bounds for $j(w)$}

By \eqref{Jn}, we have 
\begin{equation}	\label{enc1}
	\dfrac{\Re(J(w_{r_{k-1}-1}))}{q_{r_{k-1}-1}} \leq 	\dfrac{\alpha J(w_{r_{k-2}-1}) + \beta J(w_{r_{k-1}-1})}{\alpha q_{r_{k-2}-1}  +  \beta q_{r_{k-1}-1}} \leq \dfrac{\Re(J(w_{r_{k-2}-1}))}{q_{r_{k-2}-1}}
\end{equation}
if $\frac{\Re(J(w_{r_{k-2}-1}))}{q_{r_{k-2}-1}}  \geq \frac{\Re(J(w_{r_{k-1}-1}))}{q_{r_{k-1}-1}}$,
 and otherwise
the inequalities are reversed in \eqref{enc1}.
We have similar bounds for the imaginary parts.


Let us choose $k_0=12$.
We  checked computationally that for all $\ell\leq 12$,
\begin{align}	
1251.36168 &\leq \dfrac{\Re(J(w_\ell))}{q_\ell}\leq 1359.5674,\label{compreal}\\
-0.4813 &\leq \dfrac{\Im(J(w_\ell))}{q_\ell}\leq 0.	\label{compimag}
\end{align}
By \eqref{enc1},
\begin{align*}
1251.36168 + \dfrac{\Re(\delta)}{q_{n}} &\leq \dfrac{\Re(J(w_{n}))}{q_{n}}\leq 1359.5674 + \dfrac{\Re(\delta)}{q_{n}}\\
-0.4813 + \dfrac{\Im(\delta)}{q_{n}} &\leq \dfrac{\Im(J(w_{n}))}{q_{n}}\leq  \dfrac{\Im(\delta)}{q_{n}}.
\end{align*}
From the bounds \eqref{erdeltak}, \eqref{erdeltaj}, \eqref{erdeltan}, \eqref{deltam1} with $k_0=12$ we obtain
\begin{align}
	\dfrac{|\Re(\delta)|}{q_n} &\leq 1.41173,	\label{deltafinalreal}\\
	\dfrac{|\Im(\delta)|}{q_n} &\leq 1.23611. 	\label{deltafinalimag}
\end{align}
By  \eqref{asimqn}, \eqref{bound delta1} and \eqref{bound delta2} we have
\begin{align}	
3206.24623 \leq \dfrac{\Re(J(w_n))}{\sqrt{n}}&\leq 3491.04708 \qquad (as\ n\rightarrow\infty), \label{encrealJn}\\
- 4.40533	\leq \dfrac{\Im(J(w_n))}{\sqrt{n}}&\leq 3.170734  \qquad (as\ n\rightarrow\infty) \label{encimagJn}
\end{align}
and by \eqref{loge},
\begin{align}	
681.50081 \leq \Re(j(w_n)) &\leq 742.03641\qquad (as\ n\rightarrow\infty)\label{encrealjn}\\
 -0.93637\leq \Im(j(w_n)) &\leq 0.67396\qquad (as\ n\rightarrow\infty). \label{encimagjn} 
\end{align}

\begin{remark}
	The condition $n\rightarrow\infty$ in \eqref{encrealJn}, \eqref{encimagJn}, \eqref{encrealjn} and \eqref{encimagjn} could be replaced by $n\geq n_0$ for an explicit value of $n_0$ if the asymptotic relations \eqref{asimqn} and \eqref{asimcn} are made explicit. Then one could compute the values $J(w_n)$, $j(w_n)$ for $n<n_0$ and hence obtain upper and lower bounds for all $n\geq 1$.
	\end{remark}
 
 \section{Appendix}
 
The first  table below shows the values of $J(w(p/q)/q)$ and $j(w(p/q))$ for the first 40  Farey fractions $p/q$ (ordered as real numbers) among all   Farey fractions up to the level $2^{12}=4096$ on the Farey tree. The second table shows the last $40$ values. The programs used for the computations were done in collaboration with Don Zagier.

	 \begin{longtable}[c]{ |c||c|c| } 
		\hline
		$p/q$ & $J(w(p/q))/q$ & $j(w(p/q))$ \\ 
				\hline\hline
				\endfirsthead
0& 1359.56741044 & 706.324813541 \\
1/14& 1341.67984291 - 0.122490502636*I& 706.858789119 - 0.0645336432753*I\\
1/13& 1340.30387617 - 0.131912848914*I& 706.900488474 - 0.0695732206634*I\\
2/25& 1339.53333481 - 0.137189362590*I& 706.923879686 - 0.0724001664859*I\\
1/12& 1338.69858166 - 0.142905585739*I& 706.949252302 - 0.0754665750541*I\\
3/35& 1338.10232944 - 0.146988601101*I& 706.967396097 - 0.0776593436020*I\\
2/23& 1337.79124132 - 0.149118869986*I& 706.976869216 - 0.0788042174199*I\\
3/34& 1337.47100355 - 0.151311793837*I& 706.986625824 - 0.0799833523780*I\\
1/11& 1336.80141549 - 0.155896998255*I& 707.007041981 - 0.0824507472184*I\\
4/43& 1336.27197388 - 0.159522502181*I& 707.023200307 - 0.0844035587194*I\\
3/32& 1336.08997833 - 0.160768769155*I& 707.028757860 - 0.0850752157429*I\\
5/53& 1335.94232156 - 0.161779891418*I& 707.033267994 - 0.0856202872524*I\\
2/21& 1335.71732077 - 0.163320649151*I& 707.040142609 - 0.0864511175163*I\\
5/52& 1335.48799304 - 0.164891036840*I& 707.047151949 - 0.0872982300007*I\\
3/31& 1335.33264200 - 0.165954847856*I& 707.051901658 - 0.0878722552532*I\\
4/41& 1335.13561141 - 0.167304071583*I& 707.057927362 - 0.0886004905281*I\\
1/10& 1334.52481658 - 0.171486665136*I& 707.076619003 - 0.0908594653933*I\\
5/49& 1334.01374406 - 0.174986346915*I& 707.092272845 - 0.0927512885887*I\\
4/39& 1333.88269983 - 0.175883701217*I& 707.096288695 - 0.0932366184737*I\\
7/68& 1333.78827089 - 0.176530324170*I& 707.099182986 - 0.0935864040170*I\\
3/29& 1333.66128026 - 0.177399920555*I& 707.103075992 - 0.0940568878788*I\\
8/77& 1333.54913268 - 0.178167875804*I& 707.106514624 - 0.0944724589373*I\\
5/48& 1333.48137685 - 0.178631848767*I& 707.108592428 - 0.0947235689505*I\\
7/67& 1333.40350822 - 0.179165071426*I& 707.110980627 - 0.0950121914106*I\\
2/19& 1333.20678745 - 0.180512160249*I& 707.117015286 - 0.0957415017533*I\\
7/66& 1333.00708607 - 0.181879659509*I& 707.123143307 - 0.0964820951663*I\\
5/47& 1332.92635572 - 0.182432478358*I& 707.125621143 - 0.0967815506251*I\\
8/75& 1332.85531302 - 0.182918958946*I& 707.127801902 - 0.0970451031725*I\\
3/28& 1332.73606276 - 0.183735551361*I& 707.131463014 - 0.0974875617700*I\\
7/65& 1332.59846631 - 0.184677773378*I& 707.135688236 - 0.0979981950158*I\\
4/37& 1332.49433927 - 0.185390806256*I& 707.138886314 - 0.0983846943323*I\\
5/46& 1332.34720323 - 0.186398352714*I& 707.143406240 - 0.0989309436276*I\\
1/9& 1331.74231064 - 0.190540488152*I&707.161999257 - 0.101177976749*I\\
6/53& 1331.21731397 - 0.194135299692*I& 707.178150980 - 0.103129850277*I\\
5/44& 1331.10992829 - 0.194870602053*I& 707.181456402 - 0.103529297799*I\\
9/79& 1331.03788473 - 0.195363906168*I& 707.183674280 - 0.103797319799*I\\
4/35& 1330.94731569 - 0.195984059913*I& 707.186462829 - 0.104134305293*I\\
11/96& 1330.87278491 - 0.196494394766*I& 707.188757875 - 0.104411652668*I\\
7/61& 1330.83002135 - 0.196787209845*I& 707.190074826 - 0.104570801246*I\\
10/87& 1330.78283397 - 0.197110316139*I& 707.191528119 - 0.104746426028*I\\
	\hline
 \end{longtable}

\begin{center}
	\begin{tabular}{ |c||c|c| } 
		\hline
		$p/q$ & $J(w(p/q))/q$ & $j(w(p/q))$ \\ 
		\hline\hline
67/144& 1256.80214081 - 0.102528424246*I& 709.686203382 - 0.0578953566189*I\\
47/101& 1256.79136764 - 0.102325397664*I& 709.686610778 - 0.0577812408124*I\\
74/159& 1256.78161080 - 0.102141524533*I& 709.686979746 - 0.0576778888468*I
\\
27/58& 1256.76462045 - 0.101821331667*I& 709.687622275 - 0.0574979097853*I\\
61/131& 1256.74399856 - 0.101432700632*I& 709.688402163 - 0.0572794549412*I\\
34/73& 1256.72761405 - 0.101123925288*I& 709.689021820 - 0.0571058826255*I\\
41/88& 1256.70322348 - 0.100664271084*I& 709.689944295 - 0.0568474876195*I\\
7/15& 1256.58452267 - 0.0984272872898*I& 709.694434219 - 0.0555898124674*I\\
36/77& 1256.44886461 - 0.0958707343734*I& 709.699566670 - 0.0541521589596*I\\
29/62& 1256.41604412 - 0.0952522135066*I& 709.700808570 - 0.0538042897974*I\\
51/109& 1256.39285899 - 0.0948152767474*I& 709.701685917 - 0.0535585356735*I\\
22/47& 1256.36227436 - 0.0942388920864*I& 709.702843321 - 0.0532343345784*I\\
59/126& 1256.33581623 - 0.0937402736098*I& 709.703844616 - 0.0529538613138*I\\
37/79& 1256.32007532 - 0.0934436271744*I& 709.704440344 - 0.0527869914925*I\\
52/111& 1256.30220725 - 0.0931068933828*I& 709.705116596 - 0.0525975662884*I\\
15/32& 1256.25809547 - 0.0922755818348*I& 709.706786180 - 0.0521298981978*I\\
53/113& 1256.21476442 - 0.0914589837656*I& 709.708426336 - 0.0516704732988*I\\
38/81& 1256.19764598 - 0.0911363771209*I& 709.709074332 - 0.0514889627848*I\\
61/130& 1256.18276611 - 0.0908559574990*I& 709.709637605 - 0.0513311839730*I\\
23/49& 1256.15816877 - 0.0903924066955*I& 709.710568762 - 0.0510703571911*I\\
54/115& 1256.13036308 - 0.0898683927437*I& 709.711621420 - 0.0507754964001*I\\
31/66& 1256.10971946 - 0.0894793520826*I& 709.712402972 - 0.0505565755912*I\\
39/83& 1256.08111686 - 0.0889403198412*I& 709.713485890 - 0.0502532388899*I\\
8/17& 1255.97007146 - 0.0868476064334*I& 709.717690658 - 0.0490754390393*I\\
33/70& 1255.83840334 - 0.0843662462496*I& 709.722677340 - 0.0476786169351*I\\
25/53& 1255.79617017 - 0.0835703382661*I& 709.724277078 - 0.0472305134542*I\\
42/89& 1255.76295307 - 0.0829443432229*I& 709.725535379 - 0.0468780499190*I\\
17/36& 1255.71405012 - 0.0820227394093*I& 709.727388008 - 0.0463591090605*I\\
43/91& 1255.66622196 - 0.0811213906246*I& 709.729200069 - 0.0458515317763*I\\
26/55& 1255.63491625 - 0.0805314168746*I& 709.730386225 - 0.0455192770302*I\\
35/74& 1255.59641869 - 0.0798059086145*I& 709.731844962 - 0.0451106692538*I\\
9/19& 1255.48497839 - 0.0777057531246*I& 709.736068162 - 0.0439277065400*I\\
28/59& 1255.34520582 - 0.0750716597984*I& 709.741366188 - 0.0424436737723*I\\
19/40& 1255.27881384 - 0.0738204654684*I& 709.743883191 - 0.0417386346705*I\\
29/61& 1255.21459865 - 0.0726102939033*I& 709.746317939 - 0.0410566358914*I\\
10/21& 1255.09228401 - 0.0703052052080*I& 709.750956292 - 0.0397573845879*I\\
21/44& 1254.92271143 - 0.0671095140622*I& 709.757388330 - 0.0379557027688*I\\
11/23& 1254.76788430 - 0.0641917091029*I& 709.763262682 - 0.0363102351122*I\\
12/25& 1254.49538854 - 0.0590563723747*I& 709.773605298 - 0.0334131593620*I\\
1/2& 1251.36168734 & 709.892890920\\
\hline
\end{tabular}
\end{center}

\end{document}